\documentclass[11pt]{article}
\usepackage{amssymb}

\begin{document}
\title{Packing points into a unit cube in higher space.}
\author{\'A.G.Horv\'ath\\ Department of Geometry, \\
Budapest University of Technology and Economics,\\
H-1521 Budapest,\\
Hungary}
\date{September 15, 2009}

\maketitle

\begin{abstract}
In this paper using the concept of the extended Hamming code we give a construction for dense packing of points at distance at least one in such unit cubes which dimension are a power of two.
\end{abstract}



\section{Introduction}The following problem was stated in \cite{moser} and later repeated in \cite{guy}, \cite{moser1}, \cite{brass}.

{\em Let $f(n)$ denote the maximum number of points that can be arranged in the n-dimensional unit cube (n-cube) so that all mutual
distances are at least 1. Obviously, $f(n) = 2^n$ for $n\leq 3$. Many have shown that $log f(n)\sim \frac{1}{2}n(log n)$. Determine the exact values of
$f(n)$ at least for small $n$.}

Any construction of a suitable point-set gives a lower bound for $f(n)$. Previously constructed sets (see \cite{bb}, \cite{bb1} and \cite{balint}, p.71) shows that $f (4)=17$, $f (5) \geq 34$, $f (6) \geq 76$, $f (7) \geq 184$, $f (8) \geq 481$, $f (9) \geq 994$,  $f (10) \geq 2452$ $f (11) \geq 5464$, $f (12) \geq 14705$. To show a good upper bound is usually much more difficult. In the paper \cite{bb2} the authors proved the following upper estimates:
$f (6) \leq 192$, $f (7) \leq 576$, $ f (8) \leq 2592$, $ f (9) \leq 11664$, $f (10) \leq 46656$, $f (11) \leq 248832$, $f (12) \leq 944784$.
More results are known for the dimension $5$, for example $f(5) \leq 44$ (see \cite{bb3}).
Better upper estimate $f(5) \leq 43$ was shown a short time ago in \cite{joos} and this was
most recently strengthened to $f(5) \leq 42$ in \cite{joos1}. The most recent result (see \cite{bb4}) says that $f(5) \leq 40$.

The best known asymptotic estimates can be found in \cite{bb1} they are $f(n) \leq {n}^{n/2}0,63901^n e^{o(n)}$ and $f(n) \geq n^{n/2}0,2419707^n O(\sqrt{n})$. The lower bound is not a constructive one.

In this paper using the extended Hamming codes we give a construction for packing
$$
\frac{3^n+2(n-1)3^{n/2}+1}{2n}
$$
points at distance at least one in the $n=2^k$-dimensional unit cube.
In the cases of $k=2,3$ the corresponding packing give the known better ones. This construction can be generalized by the Reed-Muller (briefly R-M) codes. Since the weight-distributions does not known in general, we give lower bound for the number of points of this construction. For odd $k$ it is equal to
$$
\frac{3^n+2(n-1)3^{n/2}+1}{2n}+l\cdot n^{\frac{1}{2}(\log_2n+1)}-\frac{2}{3}n[(n-1)(n-2)-3],
$$
where $l=4.768462...$, and for even $k$ is
$$
\frac{3^n+2(n-1)3^{n/2}+1}{2n}+(l-2)\cdot n^{\frac{1}{2}(\log_2n+1)}-\frac{2}{3}n[(n-1)(n-2)-3]+\frac{1}{2}\sqrt{n}.
$$

\section{Constructions using the binary extended Hamming code}

The codes that Hamming devised, the single-error-correcting binary Hamming codes and their single-error-correcting, double-error-detecting extended versions marked the beginning of coding theory. These codes remain important to this day, for theoretical and practical reasons as well as historical. In our paper we give a new application of the extended version of this code, we use it to define large point system in the unit cube at distance at least one. To the calculation of the number of points we have to use the weight distribution of this code, which is the distribution of the Hamming weights of the codewords belong to this code. So we recall to its generator function.

A simple definition of the binary Hamming code is the following one. Let $k$ be a positive integer and construct a binary $k\times (2^k-1)$ matrix $H$ such that each nonzero binary $k$-tuple occurs exactly once as a column of $H$. The codewords (of lengths $2^k-1$) of the Hamming code are the binary  $(2^k-1)$-tuples as column vectors orthogonal to the rows of $H$. It is easy to see that binary Hamming code is a linear code having minimum Hamming weight and distance $3$, meaning that every codewords contains at least three nonzero coordinates. If we add an overall parity check bit to a binary Hamming code then we have an extended Hamming code. This means that for a codeword  we glue a $2^k$-th coordinate it is zero if the weight of the codeword is even and one if it is odd, respectively. The extended binary Hamming code is a mod2 linear vertex set of the n-dimensional unit cube each of which contains at least four nonzero coordinates.

First we consider those points of the cube which are on the middle of the position vectors occurs as a codeword of the extended Hamming code. If the Hamming weight of a codeword is $j$ then the corresponding point there is $j$ coordinates with value $\frac{1}{2}$ and its other coordinates are zero. Corresponding to this codeword we collect those points which other coordinates are either zero or one. Such a way a codeword with weight $j$ generates $2^{2^k-j}$ points in the unit cube. The distances of these points to each other at least one, since the minimum weight of the code is $4$ showing that the substraction of two position vectors contains at least $4$ coordinates with absolute value $\frac{1}{2}$ or at least one with absolute value $1$. We now determine the number of points of this system.

The weight for the Hamming code are found by use of the correspondence between linear dependence relations among $r$ columns of the matrix $H$ and  codewords of weight $r$. The number of codewords of weight $j$ is denoted $W(j)$. The number of codewords that are zero linear combinations of $(j-2)$ vectors is $W(j-2)$. The number of linear combinations  consisting of $j-2$ terms that add to 0 plus one nonzero term $(n-(j-2))W(j-2)$ where $n=2^k-1$. The number of linear combinations of $(j-1)$ vectors is ${n\choose j-1}$. The number of nonzero linear combinations of $(j-1)$ vectors is ${n\choose j-1}-W(j-1)$. The number of nonzero linear combinations of $(j-1)$ vectors to which one more term may be added to form a linear dependence among $j$ vectors is therefore
$$
jW(j)={n \choose j-1}-W(j-1)-(n-(j-2))W(j-2).
$$
This recurrence relation gives a possibility to calculate the weights of this code. Of course, $W(0)=1$, $W(1)=W(2)=0$. To an explicit solution we introduce the generator function
$$
f(x)=\sum\limits_{j=0}^{n}W(j)x^j.
$$
As it can be seen in \cite{macwilliam} $f(x)$ can be determined and it is
$$
f(x)=\frac{1}{n+1}\left[(1+x)^n+n(1+x)^{(n-1)/2}(1-x)^{(n+1)/2}\right],
$$
where $n=2^k-1$.

Now the weights for the extended binary codes can be found by nothing that each codewords of odd weight has a 1 added to it, while each of even
weight has a 0 added. Thus the number of codewords of weight $j$ is 0 if $j$ is odd and $V(j):=W(j)+W(j-1)$ if $j$ is even. The odd terms and even terms can be separated by using $[f(x)+f(-x)]/2$ for the even terms and $[f(x)-f(-x)]/2$ for the odd terms. Then the generator function as we can see in \cite{peterson} is
$$
g(x)=\frac{1}{2}[f(x)+f(-x)]+\frac{x}{2}[f(x)-f(-x)]=
$$
$$
=\frac{1}{2^{k+1}}\left[(1+x)^{2^k}+ (1-x)^{2^k}+2(2^k-1)(1-x^2)^{2^{k-1}}\right].
$$
The number of our points can be calculated on the base of the numbers $V(j)$ since two points are distinct either they associated to distinct codewords or they are at least one distinct integer coordinate. Thus if we denote by $H(j)$ the number of points associated to codewords with weight $j$ we have:
$$
H(j)=2^{2^k-j}V(j).
$$
Equivalently, we have
$$
g(x)=\sum\limits_{j=0}^{2^k}V(j)x^j= \sum\limits_{j=0}^{2^k}\frac{1}{2^{2^k-j}}H(j)(x)^j
$$
$$
2^{2^k}g(x)=\sum\limits_{j=0}^{2^k}H(j)(2x)^j.
$$
Thus
$$
\sum\limits_{j=0}^{2^k}H(j)= 2^{2^k}g\left(\frac{1}{2}\right)= \frac{3^n+2(n-1)3^{n/2}+1}{2n},
$$
if $n$ denotes the dimension $2^k$ of our space.

This construction usable if $k\geq 2$ and gives point systems with cardinality 17, 481, 1351361 in the respective dimensions 4,8 and 16. The first two values give the best known results, respectively, but in the third case we can give better point system. Observe that in a 16-dimensional space certain points which all coordinates are either $\frac{1}{4}$ or $\frac{3}{4}$ can be attached to the previous point system, such that, the distances of it from the original points are at least one. To give a relatively large system with at least one pairwise distances, consider the extended Hamming code which is a linear code with minimal weight four, and adding the the points of forms $\frac{1}{4}v+\frac{3}{4}(1-v)$ to the system where $v$ is a codeword of this extended Haming code. Since the subtraction of two codewords contains at least four nonzero coordinates the subtraction of such two point there are at least four coordinates with absolute value greater or equal to $\frac{1}{2}$, showing that the distance between the points is at least one.
Thus it can be attached
$$
g(1)=\sum\limits_{j=0}^{2^4}V(j)=2^{2^4-4-1}=2048
$$
new points and we get a better point system with cardinality $1353409$.

\section{A generalization of the construction}

As we observe in the previous paragraph, to the codewords of greater weights we can correspond to points containing smaller coordinates. So respectively handing the codewords with the same weights is a natural possibility. However the difference of two codewords with the same weights can be codeword with smaller weight, it is more clearly to use linear sub-codes of the Hamming code as the set of codewords with fixed weights. Fortunately in the extended Hamming code $EH(2^{k})$ we can find nested sequences of linear codes with increasing minimal distances, $4,16,\ldots $ these are the Reed-Muller codes denoted by $RM[2^{k},1+\sum_{i=1}^r{k \choose i},2^{k-r}]$ for $r=(k-2),(k-4)\cdots ,k-2[\frac{k}{2}]$.

In our construction firstly we collect the vertices of the original cube as points. Then we determine the further points of the arrangement corresponding to the element of the R-M codes mentioned above. Let $v\in RM[2^{k},1+\sum_{i=1}^r{k \choose i},2^{k-r}]$ for an $r\in \{(k-2),(k-4),\ldots ,k-[\frac{k}{2}]\}$, the weight of $v$ is $w(v)\not =0$. To this codeword we correspond to all points of form
$$
\frac{l}{2^{\frac{k-r}{2}}}v+(1-v)\cdot \varepsilon,
$$
where $1$ is the codeword $(1,\ldots ,1)$ of length $2^k$, $(1-v)\cdot \varepsilon$ means an element of the universe code of length $2^k-w(v)$ positioning to the nonzero coordinates of the codeword $(1-v)$, $l$ is odd and $1\leq l\leq 2^{k-r-2}$.
Remark that a binary vector $v$ will be take into consideration with respect to all of the codes containing it.
The number of such points is $2^{\frac{k-r}{2}-1}\cdot2^{2^k-w(v)}$, and in every steps we are corresponding to $v$ such points in the unit cube which were not taking before. We also remark that the present construction is a simplified variation of that one which was investigated in the preceding paragraph. The reason is that the definition of the point system in such a way is more clearly.
Now the total number of points in this construction is
$$
\frac{3^n+2(n-1)3^{n/2}+1}{2n}+{\left(\sum\limits_{r=k-2[\frac{k}{2}]}^{k-4}\right)}^\ast\sum\limits_{w(v)>0\atop v\in RM[2^{k},1+\sum_{i=1}^r{k \choose i},2^{k-r}]}2^{\frac{k-r}{2}-1}\cdot2^{2^k-w(v)},
$$
where $(\sum)^\ast$ means that we take the sum only for even or only odd integers, respectively.
Unfortunately there is no general formula for the weight distribution of the Reed-Muller codes except in the cases of $r=0,1,2$ and $k-2$. These are the repetition, even weight, second order Reed-Muller and Hamming codes, respectively. Thus there is no chance to determine the exact value of this sum. An asymptotic lower bound we can get using the minimal weight elements only since their  number are known. Let $A_{2^{k-r}}$ be the number of codewords of minimal weight. Then we have:
$$
A_{2^{k-r}}=2^r\prod_{i=0}^{k-r-1}\frac{2^{k-i}-1}{2^{k-r-i}-1},
$$
and thus a lower bound of the second part of our sum is
$$
{\left(\sum\limits_{r=k-2[\frac{k}{2}]}^{k-4}\right)}^{\ast} 2^{\frac{k-r}{2}-1}\cdot2^{2^k-2^{k-r}}\cdot 2^r\prod_{i=0}^{k-r-1}\frac{2^{k-i}-1}{2^{k-r-i}-1}=
$$
$$
{\left(\sum\limits_{r=k-2[\frac{k}{2}]}^{k-4}\right)}^{\ast} 2^{2^k-2^{k-r}+\frac{k-r}{2}-1}\left(\frac{2^k}{2^{k-r}}\prod_{i=0}^{k-r-1}\frac{2^{k}-2^{i}}{2^{k-r}-2^{i}}\right).
$$
Introducing the notation $r'=k-r$ we get that this sum is:
$$
{\left(\sum\limits_{r'=4}^{2[\frac{k}{2}]}\right)}^{\ast} 2^{2^k-2^{r'}+\frac{r'}{2}-1}\left(\frac{2^k}{2^{r'}} \prod_{i=0}^{r'-1}\frac{2^{k}-2^{i}}{2^{r'}-2^{i}}\right).
$$
For odd $k$ the maximal value of $r'$ is $2[\frac{k}{2}]=k-1$ meaning that if $k\geq r'\geq 2$ then
$$
2^k-2^{r'}+\frac{r'}{2}-1\geq \frac{r'^2}{2}+\frac{r'}{2}+1,
$$
and thus the number of points is greater than
$$
\frac{3^n+2(n-1)3^{n/2}+1}{2n}
+{\left(\sum\limits_{r'=4}^{k-1}\right)}^{\ast} 2^{\frac{r'^2}{2}+\frac{r'}{2}+1} \left(\frac{2^k}{2^{r'}}\prod_{i=0}^{r'-1} \frac{2^{k}-2^{i}}{2^{r'}-2^{i}}\right)=
$$
$$
=\frac{3^n+2(n-1)3^{n/2}+1}{2n}
+{\left(\sum\limits_{r'=0}^{k-1}\right)}^{\ast} 2^{\frac{r'^2}{2}+\frac{r'}{2}+1} \left(\frac{2^k}{2^{r'}}\prod_{i=0}^{r'-1} \frac{2^{k}-2^{i}}{2^{r'}-2^{i}}\right)- \frac{2}{3}n[(n-1)(n-2)-3]=
$$
$$
=\frac{3^n+2(n-1)3^{n/2}+1}{2n}
+l\cdot n^{\frac{1}{2}(\log_2n+1)}-\frac{2}{3}n[(n-1)(n-2)-3],
$$
where $l=4.768462...$, as John Leech enumerated in \cite{leech}.

Analogously it can be seen, that the lower bound in the case when the number $k$ is even a little bit worth, it is
$$
\frac{3^n+2(n-1)3^{n/2}+1}{2n}
+l\cdot n^{\frac{1}{2}(\log_2n+1)}-\frac{2}{3}n[(n-1)(n-2)-3]- 2(n^{\frac{1}{2}(\log_2n+1)}-\frac{1}{4}\sqrt{n}),
$$
showing that the number of points for large $n$ is greater then
$$
\frac{3^n+2(n-1)3^{n/2}+1}{2n}
+(l-2)\cdot n^{\frac{1}{2}(\log_2n+1)}-\frac{2}{3}n[(n-1)(n-2)-3]+\frac{1}{2}\sqrt{n},
$$
as we stated.

\end{document}